\documentclass[12pt,leqno]{amsart}

\usepackage[utf8]{inputenc}   
\usepackage[T1]{fontenc}      
\usepackage{geometry}         
\usepackage[francais]{babel}  
\usepackage[all]{xy}
\usepackage{amssymb}

\voffset-1cm \markleft{\rm AcAa-algebras} 
\newtheorem{theorem}{Theorem}
\newtheorem{definition}[theorem]{Definition}

\newtheorem{lemma}[theorem]{Lemma}

\newtheorem{proposition}[theorem]{Proposition}

\newtheorem{defi}{D\'{e}finition}

\newcommand{\K}{\mathbb K}

\newcommand{\h}{\frak{h}}

\newcommand{\N}{\mathbb N}

\newcommand{\he}{\mathcal{H}}

\newcommand{\pf}{\noindent{\it Proof. }}

\newcommand{\ad}{\text{\rm ad}}

\usepackage{graphicx}

\title{Anti-commutative anti-associative algebras. Acaa-algebras}
\author{Elisabeth Remm}
\date{}
\address{Universit\'e de Haute-Alsace, IRIMAS UR 7499, F-68100 Mulhouse, France.}
\email{elisabeth.remm@uha.fr}

\begin{document}

\maketitle

\begin{abstract} Let $(A, \mu)$ be a nonassociative algebra over a field of characteristic zero. The polarization process allows us to associate two other algebras, and this correspondence is one-one, one commutative, the other anti-commutative. Assume that $\mu$ satisfies a quadratic  identity $\sum_{\sigma \in \Sigma_3} a_{\sigma}\mu(\mu(x_{\sigma(i)},x_{\sigma(j)}),x_{\sigma(k)}-a_{\sigma}\mu(x_{\sigma(i)},\mu(x_{\sigma(j)},x_{\sigma(k)})=0.$ Under certain conditions, the polarization of such a multiplication determines an anticommutative multiplication also verifying a quadratic identity. Now only two identities are possible, the first is the Jacobi identity which makes this anticommutative multiplication a Lie algebra and the multiplication $\mu$ is Lie admissible, the second, less classical is given by
$[[x,y],z]=[[y,z],x]=[[z,x],y].$ Such a multiplication is here called Acaa for Anticommutative and Antiassociative. We establish some properties of this type of algebras. 
\end{abstract}
\tableofcontents
\noindent{\bf Introduction}
The notion of polarization is naturally introduced in the elementary study of functions of a real variable. To such a function we uniquely associate a pair of functions formed by an even function and an odd function whose sum gives the initial function.

We easily extend this concept to the set of bilinear applications on a vector space with values in itself (or in an another space).
If $\varphi$ is such an application, the polarization process associates with it the pair $\varphi_a,\varphi_c$ defined by the identities
$$\varphi_a(X,Y)=\varphi(X,Y)-\varphi(Y,X), \ \ \varphi_c(X,Y)=\varphi(X,Y)+\varphi(Y,X).$$
The most classic example is to consider on the vector space $A$ an algebra structure defined by a multiplication $\mu$. If we assume that $A$ is an associative algebra, thet is $\mu$ satisfies the relation
$$\mu(X,\mu(Y,Z))=\mu(\mu(X,Y),Z)$$
then $\mu_c$ is a Lie bracket, that is satisfies the Jacobi identity
$$\mu(X,\mu(Y,Z))+\mu(Y,\mu(Z,X))+\mu(Z,\mu(X,Y))=0.$$
\medskip

This work is motivated by the following remark. Assume that $\varphi$ is an algebra law on $A$ (multiplication defining the structure of the algebra $A$) that is a nonassociative bilinear map satisfying an identity of the type
\begin{equation}\label{mu}
\begin{array}{ll}a_1(x_1x_2)x_3+a_2(x_2x_1)x_3+a_3(x_3x_2)x_1+a_4(x_1x_3)x_2+a_5(x_2x_3)x_1+a_6(x_3x_1)x_2 & \\
+b_1x_1(x_2x_3)+b_2x_2(x_1x_3)+b_3x_3(x_2x_1)+b_4x_1(x_3x_2)+b_5x_2(x_3x_1)+b_6x_3(x_1x_2)&=0.
\end{array}
\end{equation}
where we write $xy$ in place of $\varphi (x,y)$  Suppose that this identity reduces to an identity relative to $\varphi_a$.  For example
$$
\begin{array}{ll}(x_1x_2)x_3-(x_2x_1)x_3-(x_3x_2)x_1-(x_1x_3)x_2+(x_2x_3)x_1+(x_3x_1)x_2 & \\
-x_1(x_2x_3)+x_2(x_1x_3)+x_3(x_2x_1)+x_1(x_3x_2)+x_2(x_3x_1)+x_3(x_1x_2)&=0.
\end{array}
$$
is equivalent to write the Jacobi identity for $\varphi_a$ (We say in this case that $\varphi$ or $A$ is Lie-admissible).  The question then arises whether there is any other type of relation than the Lie admissibility which is equivalent to a relation concerning only the bracket $\varphi_a$. This problem has been addressed in \cite{ER-pol}. The following result appears: apart from the Lie admissibility there exists one (and only one) other relation, whose associated bracket verifies
\begin{equation}
\label{AaAc}
\varphi_c(x_1,\varphi_c(x_2,x_3))=\varphi_c(x_2,\varphi_c(x_3,x_1))=\varphi_c(x_3,\varphi_c(x_1,x_2))
\end{equation}
In this case we shall say that the algebra $(A,\varphi_c)$ is  an anticommutative anti-associative algebra. The antiassociativity \cite{R-Anti} is satisfied because
$$\varphi_c(x_1,\varphi_c(x_2,x_3)+\varphi_c(\varphi_c(x_1,x_2),x_3)=\varphi_c(x_1,\varphi_c(x_2,x_3)-\varphi_c(x_3,\varphi_c(x_1,x_2))=0.$$
The goal is to study these algebras that we will call Acaa-algebras.

\section{Definition and basic properties}
Let $\K$ be a field of characteristic $0$. Recall that  an algebra $(A,\mu)$ over $\K$ (often simply called an algebra) is a $\K$-vector space $A$ equipped with a bilinear product $\mu$. 
\subsection{Definition}
\begin{definition}
We say that the algebra $(A,\mu)$ is an Acaa-algebra if the product $\mu$ is anticommutatif
$$\mu(x_1,x_2)=-\mu(x_2,x_1)$$ and satisfies the quadratic identity
\begin{equation}
\label{AA}
\mu(x_1,\mu(x_2,x_3)=\mu(x_2,\mu(x_3,x_1))
\end{equation}
for any $x_1,x_2,x_3 \in A$.
\end{definition}
To clearly show the anticommutativity  of this product we will denote $\mu(x,y)$ by $[x,y]$. It is clear that we have also
$$[x_1,[x_2,x_3)]=[x_2,[x_3,x_1]]=[x_3,[x_1,x_2]].$$ 
This product is also antiassociatif \cite{R-Anti}, because
$$[x_1,[x_2,x_3)]=[x_3,[x_1,x_2]]=-[[x_1,x_2],x_3,].$$

\medskip

\noindent {\bf Examples}
\begin{enumerate}
  \item Any abelian algebra ($[x,y]=0$ for any $x,y \in A$ is an Acaa-algebra.
  \item A Lie algebra is an Acaa algebra if and only if the Lie bracket is $2$-step nilpotent.
In fact if $[,]$ is an Acaa-bracket, the Jacobi condition gives
$$[x_1,[x_2,x_3)]+[x_2,[x_3,x_1]]+[x_3,[x_1,x_2]]=3[x_1,[x_2,x_3)]=0$$
This shows that the corresponding Lie algebra is $2$-step-nilpotent.
 \item Let us denote by $(\mathcal{F}_{Acaa}(X_1,\cdots,X_n),\cdot)$ the free Acaa-algebra with $n$-generators. It is a graded algebra. Let us examine the first particular cases
$$\mathcal{F}_{Acaa}(X_1)= \K\{X\}$$
 \item Let $\mathcal{F}_{Acaa}(X_1,X_2)$ be the free Acaa-algebra with $2$ generators $X_1,X_2$.  Then
 $$\mathcal{F}_{Acaa}(X_1,X_2)=\K\{X_1,X_2\} \oplus \K\{[X_1,X_2]\}.$$
 \item Let us consider the $7$-dimensional anticommutative algebra whose bracket is given, in a basis $\{e_i,1 \leq i \leq 7\}$ by
 $$
 \left\{
 \begin{array}{l}
     [e_1,e_2]= e_4,  [e_1,e_3]= e_5,  [e_2,e_3]= e_6  \\
     \lbrack e_1,e_6]=-   \lbrack e_2,e_5]=   \lbrack e_3,e_4]= e_7.  \\
\end{array}
\right.
$$
It is a Acaa-algebra but not a Lie algebra. 
 \end{enumerate}
 
 \subsection{A characterization of the Acaa-algebras}
Assume that $(A,[,])$ is a Acaa-algebra.  By definition
$$[x_1,[x_2,x_3]]=[x_2,[x_3,x_1]]=[x_3,[x_1,x_2]].$$ 
In particular, if $x_3=x_1$, we obtain
$$[x_1,[x_2,x_1]]=[x_2,[x_1,x_1]]=0.$$ Conversely, assume that 
$$[x_1,[x_2,x_1]]=0$$
for any $x_1,x_2 \in A$. Then
$$[x_1+x_3,[x_2,x_1+x_3]]=0=[x_1,[x_2,x_3]]+[x_3,[x_2,x_1]]$$
We deduce
$$[x_1,[x_2,x_3]]=[x_3,[x_1,x_2]]$$
and $(A,[,])$ is an  Acaa-algebra.
\begin{proposition}
Let $(A,[,])$ be an anticommutative  algebra. Then $A$ is an Acaa-algebra if and only if
$$[x_1,[x_2,x_1]]=0$$
for any $x_1,x_2 \in A$. 
\end{proposition}

\section{Classification in small dimensions}
In this section we denote by $\{e_1,\cdots,e_n\}$ a basis of $A$ when $\dim A=n$.
\subsection{Dimension $2$}
We put $[e_1,e_2]=ae_1+be_2.$ Then
$$[e_1,[e_1,e_2]]=0=b[e_1,e_2], \ \ [e_2,[e_1,e_2]]=0=a[e_2,e_1]$$
this implies $[e_1,e_2]=0.$
Any $2$-dimensional Acaa-algebra is abelian.
\subsection{Dimension $3$}
\begin{lemma} Let $e$ and $f$ be two independent vectors in $A$. Then, if $[e,f]$ is not zero, then $e,f,[e,f]$ are independent.
\end{lemma}
\pf Assume that $[e,f] \neq 0$ and $[e,f]=ae+bf$. Then
$$[e,[e,f]]=0=b[e,f], \ [f,[e,f]]=0=-a[e,f]$$
this implies $a=b=0$ which contradicts the    hypothesis. 

Assume that $A$ is not abelian. There are two independent vectors, for example $e_1$ and $e_2$ such that $[e_1,e_2] \neq 0.$ From the lemma, $e_3=[e_1,e_2]$ is not in $\K\{e_1,e_2\}$ and $\{e_1,e_2,e_3\}$ is a basis of $A$. We have also
$$[e_1,e_3]=[e_1,[e_1,e_2]]=0, \ \  [e_2,e_3]=[e_2,[e_1,e_2]]=0.$$
In this case $A$ is a Lie algebra isomorphic to the $3$-dimensional Heisenberg algebra.
\begin{proposition}
Every $3$-dimensional Acaa-algebra is isomorphic to one of the algebras below
\begin{enumerate}
  \item the $3$-dimensional abelian algebra
  \item the $3$-dimensional Heisenberg algebra $\h_3$ given by $[e_1,e_2]=e_3.$
\end{enumerate}
\end{proposition}
\subsection{Dimension $4$}
We assume that $A$ is not abelian. Thus we can find a basis $\{e_1,e_2,e_3,e_4\}$ of $A$ such that $\{e_1,e_2,e_3\}$ generates a subalgebra isomorphic to the Heisenberg algebra that is $[e_1,e_2]=e_3,\ [e_i,e_3]=0$ for $i=1,2$. It remains to compute $[e_i,e_4]$ for $i=1,2,3$.  Let us put $[e_1,e_4]=a_1e_1+a_2e_2+a_3e_3+a_4e_4$. Then $[e_1,[e_1,e_4]]=0$ gives
$b_1e_3+d_1[e_1,e_4]=0$ If $d_1 \neq 0$, then $[e_1,e_4]=-(b_1/d_1)e_3$ which contradics $d_1 \neq 0$. Then $d_1=b_1=0$ and $[e_1,e_4]=a_1e_1+c_1e_2$. Now
$$[e_2,[e_1,e_4]=-a_1e_3=[e_4,[e_2,e_1]]=[e_3,e_4].$$
We deduce $a_1=0$ and $[e_1,e_4]=c_1e_3, [e_3,e_4]=0.$ A similar calculation gives $[e_2,e_4]=c_2e_3$. In conclusion
$$[e_1,e_2]-e_3, \  [e_1,e_4]=c_1e_3, \  [e_2,e_4]=c_2e_3$$
and $A$ is a Lie algebra isomorphic to the direct sum of $\he_3$ and the $1$-dimensional Lie algeba.
\begin{proposition}
Every $4$-dimensional Acaa-algebra is isomorphic to one of the algebras below
\begin{enumerate}
  \item the $4$-dimensional abelian algebra
  \item the $4$-dimensional $2$-step-nilpotent Lie algebra  $\he_3 \oplus \K$ given by $[e_1,e_2]=e_3.$
\end{enumerate}
\end{proposition}
\subsection{Dimension 5}
As in dimension $4$, we assume that $A$ is not abelian. Thus we can find a basis $\{e_1,e_2,e_3,e_4,e_5\}$ of $A$ such that $\{e_1,e_2,e_3\}$ generates a subalgebra isomorphic to the Heisenberg algebra that is $[e_1,e_2]=e_3,\ [e_i,e_3]=0$ for $i=1,2$. In a first time we assume that $[e_1,e_4]=e_5$. This implies $[e_1,e_5]=[e_4,e_5]=0$ and 
$$[e_2,e_5]=[e_3,e_4]=[e_1,[e_4,e_2]]=ae_3+be_5.$$ 
We deduce $[e_2,[e_3,e_4]]=b[e_2,e_5]=[e_4,[e_2,e_3]]=0$ and $b=0$ or $[e_2,e_5]=0$. If $b=0$, then $[e_3,e_4]=ae_3$ this giving $a=0$. In any case, $[e_2,e_5]=[e_3,e_4]=0$ and $A$ is isomorphic to the $2$-step nilpotent Lie algebra given by
$$[e_1,e_2]=e_3, [e_1,e_4]=e_5, [e_2,e_4]=c e_3 +d e_5$$
by isomorphism we can take $c=d=0$. If the hypothesis $[e_1,e_4]=e_5$ is not satisfied, that is $[e_1,e_4] \in \K\{e_1,e_2,e_3,e_4\}$, we can consider that $A_1=\K\{e_1,e_2,e_3,e_4\}$ is subalgebra of $A$. Then it is isomorphic to $\he_3 \oplus \K$. In this case, we have to compute $[[e_i,e_5]$ and $A$ , if it is not of the previous case,  is isomorphic to $\he_5$.

\begin{proposition}
Every $5$-dimensional Acaa-algebra is isomorphic to one of the algebras below
\begin{enumerate}
  \item the $5$-dimensional abelian algebra
  \item the $5$-dimensional $2$-step-nilpotent Lie algebra  $\he_3 \oplus \K^2$ given by $[e_1,e_2]=e_3.$
  \item the $5$-dimensional $2$-step-nilpotent Lie algebra defined by
  $$[e_1,e_2]=e_3, \ [e_1,e_4]=e_5.$$
  \item the $5$-dimensional $2$-step-nilpotent Lie algebra $\he_5$.
\end{enumerate}
\end{proposition}
In this case also, any Acaa-algebra is a Lie algebra. In the next section, we shall see a $7$-dimensional Acaa-algbera which is not a Lie algebra. 
\subsection{The free Acaa-algebras}
In the first section we have given, as example, the structure of the free Acaa-algebras with $1$ or $2$-generators. Let us examine the general case. We consider $n$ generators $X_1,\cdots, X_n$ and we denote by
$$\mathcal{F}_{Acaa}(X_1,\cdots, X_n)=\oplus_{k \geq 1}\mathcal{F}^k_{Acaa}(X_1,\cdots, X_n)$$
the free Acaa-algebras with $n$ generators, where $\mathcal{F}^k_{Acaa}(X_1,\cdots, X_n)$ denotes the homogeneous component of degree $k$. We have
\begin{enumerate}
  \item $\mathcal{F}^1_{Acaa}(X_1,\cdots, X_n)=\K\{X_1,\cdots, X_n\}$
  \item $\mathcal{F}^2_{Acaa}(X_1,\cdots, X_n)=\K\{X_{ij}, \ 1 \leq i < j \leq n\}$ where $X_{ij}=[X_i,X_j]$,
  \item $\mathcal{F}^3_{Acaa}(X_1,\cdots, X_n)=\K\{X_{ijk}, \ 1 \leq i < j < k \leq n\}$ where $X_{ijk}=[X_i,[X_j,X_k]]$,
  \item $\mathcal{F}^k_{Acaa}(X_1,\cdots, X_n)=0$ for any $k \geq 4$.
\end{enumerate}
\pf In fact since $[X_i,[X_j,X_k]]=[X_j,[X_k,X_i]]=[X_k,[X_i,X_j]]$ the components of degree $3$ are linear combination of the vectors $X_{ijk}$ with $1 \leq i< j< k \leq n.$ Concerning the degree $4$, from \cite{R-Anti}, any anti-associative algebra is $4$-step nilpotent and any product of degree $4$ is equal to $0$.

\medskip

\noindent{Remark} We have seen that $\mathcal{F}_{Acaa}(X_1, X_2)$ is a $3$-dimensional Lie algebra isomorphic to $\he_3$. For $n=3$ the generators are $\{X_1,X_2,X_3,X_{12},X_{13},X_{23},X_{123}\}$ and we have
$$\left\{
\begin{array}{l}
    \lbrack X_1,X_2]=X_{12},  \lbrack X_1,X_3]=X_{13},   \lbrack X_2,X_3]=X_{23},    \\
     \lbrack X_1,X_{23}] = -\lbrack X_2,X_{13}]=\lbrack X_3,X_{12}]=X_{123}.
\end{array}
\right.
$$
This algebra is not $2$-step nilpotent and it has no  Lie algebra  structure, for this bracket.
\section{Representations of Acaa-algebras}

\subsection{The operator $\ad X$}
Let $A$ be a Acaa-algebra. Let $\ad X$ be the linear map on $A$ defined by $\ad X(Y)=[X,Y].$ Since $A$ is an Acaa-algebra, , for any $X,Y,Z \in A$, we have
$$[X,[Y,Z]]=[Y,[[Z,X]]=[Z,[X,Y]]$$
that is
$$[[X,Y],Z]=-[X,[Y,Z]]=[Y,[X,Z]].$$
This is equivalent to
$$\ad [X,Y]=-\ad X \circ \ad Y = \ad Y \circ \ad X.$$
We deduce
\begin{lemma}
The operators $\ad X$ on $A$ satisfy
\begin{enumerate}
  \item $\ad X \circ \ad Y + \ad Y \circ \ad X =0$
  \item $2\ad [X,Y]=-\ad X \circ \ad Y + \ad Y \circ \ad X= - [\ad X, \ad Y]_{gl(n)}$
\end{enumerate}
where $[,]_{gl(n)}$ denotes the classical bracket on the space of linear operators.
\end{lemma}

\medskip

\noindent{\bf Applications} The linear map
$$\ad : A \rightarrow End(A)$$
is a weighted linear representation of $A$. It satisfies
$$\ad \lbrack X,Y]=-\frac{1}{2}\lbrack \ad X,\ad Y].$$
Since $A$ is also characterized by the relation
$$[X,[X,Y]]=0$$
then $(\ad X)^2=0$ for all $X \in A.$ If we consider $\ad$ as a matricial representation of $A$, the matrices which represent the elements of $A$ are square null matrices. 

\medskip

Along the same lines, we have
$$\ad X [Y,Z]=-[Y,\ad X(Z)]=-[\ad X (Y), Z].$$
So,
$$2 \ad (X) [Y,Z]=-[Y, \ad X(Z)]-[\ad X(Y),Z].$$
In \cite{R-Anti}, we have introduced the notion of anti-derivation: Let $A$ be an algebra. A linear map $f \in End(A)$ is an anti-derivation if this map satisfies
$$f(xy)=-xf(y)-f(x)y=0$$
for any $x,y \in A.$
 As for the representations, we are here led to put coefficients on such a relation:
\begin{definition}
Let $B$ be an algebra. A linear map $f \in End(B)$ is an weighted anti-derivation if there exists $k \in \N^*$ such that
$$kf(xy)=-xf(y)-f(x)y=0$$
for any $x,y \in A.$
\end{definition}
In these conditions
\begin{lemma}
Let $A$ be an Acaa algebra. For any $X \in A$, the linear map $\ad X$ is a weighted anti-derivation of weight $2$. 
\end{lemma}

\medskip

\noindent{\bf Remark.} In \cite{R-Anti}, we have seen that in any anti-associative algebra, the operator $L_X+R_X$ is an anti-derivation. But in case of anti-commutative anti-associative algebra we have  $L_X+R_X=0$.

\subsection{Representations of Acaa algebras}
Since an Acaa algebra is a nonassociative algebra, the notion of representation  must be taken in the framework of the associative representations of non-associative algebras introduced in \cite{S}. Take for example the map $\ad$. It satisfies
$$\ad [X,Y]=-\ad X \circ \ad Y=\ad Y \circ \ad X.$$
So $\ad: A \rightarrow End(A)$ is a linear mapping which gives a representation of $A$ in the associative algebra $End(A)$.

\begin{definition}
Let $A$ be an Acaa algebra and let $End(A)$ be the space of linear endomorphisms of the vector space $A$. A representation of $A$ is a linear map
$$\rho : A \rightarrow End(A)$$ which satisfies the identity
$$\rho[X,Y]=-\rho(X)\circ \rho(Y)=\rho(Y) \circ \rho(X)$$
for any $X,Y \in A$.
\end{definition}
The representation $\rho$ is said to be faithful if it is injective. Of course, the linear map $\ad$ is an associative representation of $A$ which is not faithful. 

Let us note that, if $\rho$ is a representation of $A$, then $\rho(X)^2=0.$  Then, identifying $End(A)$ with the space of square matrices of order $n=\dim A$, $\rho(X)$ is considered as a matrix $\widetilde{X}$ satisfying $\widetilde{X}^2=0.$

\medskip

\noindent{\bf Example: $\dim 3$. The Heisenberg algebra}
We have seen that any $3$-dimensional Acaa algebra is isomorphic to the Heisenberg Lie algebra $\h_3$. It is defined by $[e_1,e_2]=e_3$ and it is isomorphic to the free Acaa algebra with $2$ generators.  For any $X \in A$, the matrix $\widetilde{X}=\rho (X)$ is $2$-nilpotent. Thee is a nonsingular matrix $P$ such that 
$$P^{-1} \widetilde{X}P=N_1=\begin{pmatrix}
     0 & 0 & 0  \\
    1  & 0 & 0 \\
    0 & 0 & 0 
\end{pmatrix}$$
Let  $\widetilde{Y}$ be a matrix which satisfies $\widetilde{X}\widetilde{Y}=-\widetilde{Y\widetilde{X}}$.  Then 
$$P N_1 P^{-1} \widetilde{Y}=-\widetilde{Y} P N_1 P^{-1}$$
and
$$N_1 P^{-1} \widetilde{Y}P=-P^{-1}\widetilde{Y} P N_1$$
and this shows that $Y_1=P^{-1} \widetilde{Y}P$ satisfies $N_1Y_1=-Y_1N_1.$. If we compute $Y_1$ we find
$$Y_1=\begin{pmatrix}
     a_{11} & 0 & 0   \\
     a_{21} & -a_{11} & a_{23}\\
     a_{31} & 0 & a_{33} 
\end{pmatrix}$$
Since $Y_1^2=0$, we obtain $a_{11}=a_{33}=a_{23}a_{31}=0$ and
$$Y_1=\begin{pmatrix}
    0 & 0 & 0   \\
     a_{21} & 0 & 0\\
     a_{31} & 0 & 0
\end{pmatrix} \ \ {\rm or} \ \  \begin{pmatrix}
    0 & 0 & 0   \\
     a_{21} & 0 & a_{23}\\
     0 & 0 & 0
\end{pmatrix}$$
In both cases we have $[\widetilde{X},\widetilde{Y_1}]=0$ and this associative representation is not faithful. We deduce that any associative representation of the Acaa algebra $\h_3$ is not faithful. This is not the case for the Lie representations of this Lie algebra.
\medskip

\noindent{\bf Example: $\dim 7$. The  free Acaa algebra with $3$ generators}
This Acaa algebra is defined in a graded basis $\{e_1,e_2,e_3,e_{12},e_{13},e_{23},e_{123}\}$ by
$$
\left\{
\begin{array}{l}
     [e_1,e_2]=e_{12}, \   [e_1,e_3]=e_{13}, \  [e_2,e_3]=e_{23}  \\
    \lbrack e_1,e_{23}]=-  \lbrack e_2,e_{13}]= \lbrack e_3,e_{12}]=e_{13}\\
\end{array}
\right.
$$
If $X=\sum x_ie_i$, then $\d X$ is represented by the matrix
$$\widetilde{X}=
\begin{pmatrix}
    0  & 0 & 0 & 0 & 0 & 0 & 0   \\
     0  & 0 & 0 & 0 & 0 & 0 & 0   \\
      0  & 0 & 0 & 0 & 0 & 0 & 0   \\
       -x_2 & x_1 & 0 & 0 & 0 & 0 & 0   \\ 
        -x_3  & 0 & -x_1 & 0 & 0 & 0 & 0   \\
         0  & -x_3 & x_2 & 0 & 0 & 0 & 0   \\
         -x_6  & x_5 & -x_4 & x_3 & -x_2 & x_1 & 0   \\
\end{pmatrix}
$$

\subsection{Acaa Admissible algebra}
A $\K$-algebra $(B, \cdot)$ is called Acaa-admissible if the bracket on $B$ defined by $[b_1,b_2]=b_1 \cdot b_2 -b_2 \cdot b_1$ is an Acaa-bracket. 
\begin{proposition}
The algebra $(B,\cdot)$ is Acaa-admissible if and only if the product $b_1 \cdot b_2$ satisfies the identity
$$\rho_B(b_1,b_2,b_3)-\rho_B(b_2,b_1,b_3)+\rho_B(b_1,b_3,b_2)-\rho_B(b_3,b_1,b_2)=0$$
where $\rho_B(b_1,b_2,b_3)=(b_1\cdot b_2) \cdot b_3 - b_3 \cdot (b_1 \cdot b_2)$.
for any $b_1,b_2,b_3 \in B$.
\end{proposition}
\pf See \cite{R-Anti}.
 
 \medskip
 
 In particular any algebra $(B, \cdot)$ which satisfies the identity
 $$\rho_B(b_1,b_2,b_3)=0$$
 is Acaa-admissible.  We will call  such an algebra a $\rho$-associative algebra.
 \begin{proposition}
 Let $B$ be a $\rho$-associative algebra. Then the bracket on $B$ given by $[x,y]=xy-yx$ provides  $B$ with a Acaa-algebra structure which is also a $2$-step nilpotent Lie algebra.
 \end{proposition}
 \pf In fact
 $$[[x,y],z]=(xy)z-(yx)z-z(xy)+z(yx)=\rho(x,y,z)-\rho(y,x,z).$$
 
 Now, let $F=\{X_1,\cdots,X_n\}$ be a finite set of cardinality $n$ and let $\widetilde{\K}(F)$ the nonassociative free algebra generated by $F$.  Let us consider now the quotient space
 $$\widetilde{T}(F)=\widetilde{\K}(F) / \sim$$
 where $\sim$ is the equivalence relation given by $(X_iX_j)X_k-X_k(X_iX_j)$.  Then $\widetilde{T}(F)$ is a $\rho$-associative algebra. As consequence $\widetilde{T}(F)$ is also an infinite dimensional $2$-step nilpotent Lie algebra. For example, if $n=2$, this Lie algebra is isomorphic to $\h_3 \oplus Ab_\infty$ where $Ab-\infty$ is an infinite dimensional abelian Lie algebra. If $n=3$, let us consider the $6$-dimensional $2$-step nilpotent Lie algebra $n_6$ given by
 $$[e_1,e_2]=e_4,\ [e_1,e_3]=e_5, \ [e_2,e_3]=e_6$$
 then  $\widetilde{T}(F)=n_6 \oplus Ab_\infty$. The non abelian part corresponds to the free 2-step nilpotent Lie algebra with $3$ generators.  This is also true for any $n \geq 3$.
 
 \medskip

  Let $A$ be an Acaa algebra and let us consider the $\rho$-associative algebra $\widetilde{T}(A)$ (here $F$ means a given basis of $A$).  Let $ \widetilde{U}(A)$ be the quotient ring of $\widetilde{T}(A)$ by the ideal generated by elements of the form
 $$[X_i,X_j]-X_iX_j +X_jX_i$$
 where $[,]$ is the anticommutative multiplication of $A$. 
Then $ \widetilde{U}(A)$  is also a $\rho$-associative algebra.  It can be also consider as an Acaa algebra which is also a $2$-step nilpotent Lie algebra.

\medskip

\noindent{\bf Examples}
\begin{enumerate}
  \item $\dim A =1$. Then $\widetilde{T}(A)$ is the graded space given by
  $$\widetilde{T}^1(A)=\K\{X\}, \ \widetilde{T}^2(A)=\K\{X^2\}, \ \widetilde{T}^3(A)=\K\{X^3\}, \ \widetilde{T}^4(A)=\K\{X^4=X^3X, X^2X^2\}$$
  and   $ \widetilde{U}(F)$  is abelian. 
  \item $\dim A =2$. Then $\widetilde{T}(A)=\widetilde{T}(X_1,X_2)$ is the graded space given by
  $$\widetilde{T}^1(A)=\K\{X_1,X_2\}, \ \widetilde{T}^2(A)=\K\{X_1,X_2, X_1^2,X_2^2,X_1X_2,X_2X_1\}, \cdots$$
  and if in  $ \widetilde{U}(A)$ we put $[X_1,X_2]=X_1X_2-X_2X_1=X_{12}$, then $[X_i ,X_{12}]=0$  and $\widetilde{U}(A)$ is  isomorphic to the free Acaa algebra with $2$-generators.
  \end{enumerate}


\section{Cohomology and operads}

\subsection{Cohomology of Acaa algebras}

Let $A$ be an Acaa algebra.  We saw previously that all products of order 4 were zero. We have a natural filtration
$$\{0\} \subset A^3 \subset A^2\subset A^1=A$$
where $A^i$ is the subspace of $A$ generated by the product of degree greater than $3$ of elements of $A$. Let us put $A_i= A^i-A^{i-1}$.  This space contains only the product of degree $i$. 

For any $X \in A_i$, we consider the linear map $g_X$ defined by
$$g_X(Y)=(-1)^{i+j}ij[X,Y], \ \ \forall Y \in A_i.$$
In particular $g_X(Y)=0$ as soon as $i+j \geq 4$.
Now, for every $f \in End(A)$, we put
$$\delta^1(f)(u,v)=f[u,v]-[u,f(v)]-[f(u),v]$$
for any $u,v \in A$.
In particular we have, for any $X \in A_i$, $ Y \in A_1$ and $Z \in A_1$:
$$
\begin{array}{ll}
    \delta^1(g_X)(Y,Z) &=g_X[Y,Z]-[Y,g_X(Z)]-[g_X(Y),Z]\\
      &=  -2i[X,[Y,Z]]-i[Y,[X,Z]]-i[[X,Y],Z]\\
      & =-i(2[X,[Y,Z]]-[Y,[Z,X]]-[Z,[X,Y]]\\
      &=-i([ X,[Y,Z]]-[Y,[Z,X]])-i([ X,[Y,Z]]-[Z,[X,Y]])=0.
\end{array}
$$
Now, if $\varphi$ is a skew-symmetric bilinear map on $A$ with values in $A$, we consider the trilinear map
$$\delta^2(\varphi)(X,Y,Z)=\varphi(X,[Y,Z])+[X,\varphi(Y,Z)]-\varphi(Y,[Z,X])-[Y,\varphi(Z,X)].$$
In particular we have 
$$\delta^2 \circ \delta^1 =0.$$
This trilinear map also satisfies
$$\delta^2(\varphi)(X,Y,Z)+\delta^2(\varphi)(Y,Z,X)+\delta^2(\varphi)(Z,X,Y)
=0.$$
Thus, if we want to construct something analogous to Chevalley Eilenberg complex for Lie algebras, we can define the first cochain spaces as follows:
$$\mathcal{C}^0(A,A)=A, \ \mathcal{C}^1(A,A)=End A, \ \mathcal{C}^2(A,A)=\{\varphi:A \times A \rightarrow A, {\rm bilinear \ skew-symmetric}\}$$
and 
$\mathcal{C}^3(A,A)$ the space of trilinear maps $\Psi : A \times A \times A \rightarrow A $ satifying the following identity
$$  \psi(X,Y,Z)=\psi(Y,X,Z).$$
with
$$\mathcal{C}^0(A,A)\rightarrow^{\delta^0} \mathcal{C}^1(A,A)$$
$\ \ \ \ \ \ \ \ \ \ $
\xymatrix{
    \mathcal{C}^0(A,A) \ar[r]^{\delta^0}  & \mathcal{C}^1(A,A) \ar[r]^{\delta^1}& \mathcal{C}^2(A,A) \ar[r]^{\delta^2}& \mathcal{C}^3(A,A) \ar[r]^{\delta^3} & \cdots
 &&&& }
 
\noindent and 
$$
\left\{
\begin{array}{ll}\
\delta^3(\psi)(X_1,X_2,X_3)= &\psi(X_1,X_2,[X_3,X_4]) + \psi(X_1,[X_3,X_4],X_2)+\psi(X_2,[X_3,X_4],X_1)\\
&[X_1,\psi(X_2,X_3,X_4)]+[X_1,(\psi(X_2,X_4,X_3)]+[X_1,(\psi(X_4,X_3,X_2)]
\end{array}
\right.
$$

\subsection{Quadratic operad associated with Acaa algebras}

Let $\mathcal{A}caa=\oplus_{n \geq 1} \mathcal{A}caa(n)$ be the quadratic operad associated with  Acaa algebras. We have\begin{enumerate}
  \item $\mathcal{A}caa(1)=\K$
  \item $\mathcal{A}caa(2)=\K\{x_1x_2\}$ because we have the anticommutativity.
  \item $\mathcal{A}caa(3)=\K\{x_1(x_2x_3)\}$. in fact 
  $$x_1(x_2x_3)=x_2(x_3x_1)=x_3(x_1x_2)=-x_1(x_3x_2)=-x_2(x_1x_3)=-x_3(x_2x_1)=-(x_2x_3)x_1 \cdots $$
    \item $\mathcal{A}caa(n)=\{0\}$ for $n \geq 4$ because all the products of degree $4$ are null. 
\end{enumerate}
The Poincar\'e or generating series of $\mathcal{A}caa$ is the series
$$g_{\mathcal{A}caa}(t)=\sum \frac{(-1)^n}{n!}\dim (\mathcal{AA}ss (n))t^n=-t+\frac{1}{2}t^2-\frac{t^3}{6}.$$
Since $\mathcal{A}caa(1)=\K$, $\mathcal{A}caa$ admits a minimal model, unique up to isomorphism. The generating series $g_M$ associated with  this minimal model is the formal inverse of $g_{\mathcal{A}caa(t)}$ taken with the opposite sign
$$g_{\mathcal{A}caa}(-g_M(t))=t.$$
We deduce
$$g_M(t)=-t + \frac{t^2}{2} - \frac{t^3}{3}+\frac{5t^4}{24} -\frac{t^5}{12}-\frac{7t^6}{144}+\cdots$$
and this series cannot be the Poincare Series of the dual operad. In fact, if $\mathcal{R}$ is the ideal of relations which defines 
$\mathcal{A}caa$ and which is generated as a $\Sigma_3$-module by the relation $(x_1x_2)x_3-(x_2x_3)x_1$, the ideal $\mathcal{R}^\bot$ which generates the dual operad $\mathcal{A}caa^{!}$ is the orthogonal of $\mathcal{R}$ for the classical inner product
\begin{eqnarray}
\label{pairing}
\left\{
\begin{array}{l}
<(x_i \cdot x_j)\cdot x_k,(x_{i'} \cdot x_{j'})\cdot x_{k'}>=0, \ {\rm if} \ (i,j,k) \neq (i',j',k'), \\
<(x_i \cdot x_j)\cdot x_k,(x_i \cdot x_j)\cdot x_k>={\varepsilon(\sigma)},  
\qquad   {\rm with} \ \sigma =
\left(
\begin{array}{lll}
1 & 2 & 3\\
i &j &k 
\end{array}
\right)
\\
<x_i \cdot (x_j\cdot x_k),x_{i'} \cdot (x_{j'}\cdot x_{k'})>=0, \ {\rm if} \ ( i,j,k) \neq ( i',j',k'), \\
<x_i \cdot (x_j\cdot x_k),x_i \cdot (x_j\cdot x_k)>=-{\varepsilon(\sigma)} 
\qquad {\rm with} \ \sigma =
\left(
\begin{array}{lll}
1 & 2 & 3\\
i &j &k 
\end{array}
\right)
, \\
<(x_i \cdot x_j)\cdot x_k,x_{i'} \cdot (x_{j'}\cdot x_{k'})>=0,
\end{array}
\right.
\end{eqnarray}
where $\varepsilon(\sigma)$ is the signature of $\sigma$. Then i$\mathcal{R}^\bot$ is generated as a $\Sigma_3$-module by the relation $(x_1x_2)x_3+(x_2x_3)x_1$. But this relation is equivalent to $(x_1x_2)x_3=0.$ In fact In fact any algebra associated with this operad $\mathcal{A}caa^{!}$ have to satisfy
$$(x_1x_2)x_3+(x_2x_3)x_1=0, \ (x_2x_3)x_1+(x_3x_1)x_2=0, \ (x_3x_1)x_2+(x_1x_2)x_3=0$$
and $(x_ix_j)x_k=0$ for any $i,j,k$. 
Then any $\mathcal{A}caa^{!}$-algebra is a $2$-step nilpotent Lie algebra and  the generating series of $\mathcal{A}caa^{!}$ is 
$$-t+\frac{t^2}{2}.$$ We deduce
\begin{theorem} The quadratic operad $\mathcal{A}caa$ is not a Koszul operad.
\end{theorem}

\end{document}